\documentclass[12pt]{amsart}
\usepackage{amssymb}
\usepackage{amsthm}
\usepackage{amsfonts}
\usepackage{url}
\usepackage{tikz}
\newtheorem{lemma}{\bf Lemma}[section]
\newtheorem{proposition}[lemma]{\bf Proposition}
\newtheorem{fact}[lemma]{\bf Fact}

\begin{document}
\def\Sym{\text{Sym}(\mathbb{N})}
\def\N{\mathbb{N}}
\def\PN{\mathcal{P}(\N)}
\def\str{\text{stretch}}
\parskip = 2mm
\parindent = 0mm
\title{Knuth's non-associative ``group'' on $\PN$}
\begin{abstract}Donald Knuth introduced in \cite{taocp4a} a fast approximation 
to the addition of integers (given in binary)
in terms of bit-wise operations by 
	$$a + b \; \approx \; a \oplus b \oplus ((a\land b) \ll 1).$$
Generalizing this to infinite bit-strings
we get a binary operation on $\PN$, the power-set of $\N$ (which we identify 
with the collection of infinite bit-strings). We show that this operation 
is ``group-like''
in that it has a neutral element, inverses, but it is not associative. 
There are a lot of questions left, which the author has not been able to 
answer.
\end{abstract}
\author{Dominic van der Zypen}
\maketitle
\section{Introduction}
Addition of integers  is an important operation in computer science (and in 
daily life).
Knuth \cite{taocp4a} noted that for integers $a,b$ given in binary, we have 
$$a + b = (a\oplus b) + ((a \land b) \ll 1),$$
where $\oplus$ denotes bit-wise {\tt XOR}, $\land$ is bit-wise ${\tt AND}$ and 
$\ll 1$ means shifting to left by $1$ position. 

This identity can be used for an approximation of $+$ using exclusively
bit-wise operations\footnote{These
operations are very fast operations in computers, often using only $1$ or a very 
low number of CPU-cycles}:
	$$a + b \; \approx \; a \oplus b \oplus ((a\land b) \ll 1).$$

Note that $((a\land b)\ll 1)$ is used to simulate the {\em carry-bit 
propagation}.

This approximation is not only of academic interest;
it is used in the cryptographic scheme NORX \cite{norx},
for instance.
\section{The binary operation $\oplus$ on $\PN$}
Let $\N = \{0,1,2,3,\ldots\}$ be the collection of non-negative integers 
and $\PN$ be
the power-set of $\N$, that is the collection of all subsets of $\N$.
By slight abuse of notation, we are going to define an operation 
$\oplus:\PN\times\PN \to \PN$ and will not use $\oplus$ any more 
as bit-wise {\tt XOR} on finite bit-strings.

For any set $A\in\PN$, let $A+1 = \{a+1: a\in A\}$, so $A+1$ simulates the 
{\em left-shift}. Moreover, given $A, B\in\PN$, we let 
$$A\, \triangle \, B = (A\setminus B)\cup (B\setminus A)$$ be the symmetric difference
of $A, B$. Note that $A\, \triangle \, B$ plays the role of bit-wise {\tt XOR}.

Finally, we define for all $A,B \in \PN$: $$A\oplus B := (A \, \triangle \,
B)\, \triangle \, ((A \cap B) + 1).$$
\section{Basic properties of $\oplus$}

\subsection{Commutativity} The definition is clearly symmetric on the two 
variables,
so $\oplus$ is commutative.

\subsection{Neutral element} It is easy to see that the empty set $\emptyset 
\in \PN$ is the neutral element with respect to $\oplus$.

\subsection{Non-associativity} Let $A = B = \{0\}$ and $C = \{1\}$. Then
$(A \oplus A)\oplus C = \{2\}$, but $A \oplus (A \oplus C) = \emptyset$. 

\section{Inverse elements in $\oplus$}

The goal of this section is to show that for every $A\in \PN$ there is $A'\in 
\PN$ such
that $A \oplus A'  = (A\oplus A') \oplus ((A\cap A') + 1) = \emptyset$.

First, a basic observation will be useful later:
\begin{fact}\label{basicfact}
{\em For any sets $X, Y$ we have $X\, \triangle \, Y = \emptyset$ if and 
only if $X = Y$, so $A \oplus A' = \emptyset$ amounts to saying $A\, 
\triangle \, A' =  (A\cap A')+1$.}
\end{fact}

Let us first consider a few examples:
\begin{itemize}
	\item Let $A = \{0\} \in \PN$. Then let $A' = \{0,1\}$. 
	\item More generally, let $A = \{n\}$ for some $n\in\N$. Then $A' = 
	\{n,n+1\}$.
	\item Let $A = \{3, 4, 5\}$. Then $A' = \{3, 5, 6\}$. 
\end{itemize}
Note that always we need $\min(A)\in A'$ for $A \neq \emptyset$. Now we are 
ready to construct $A'$ for general $A\in \PN$. 

We assume that $A\in\PN\setminus \{\emptyset\}$ for the remainder
of this section.

First, for $a\leq b\in\N$ we let $[a,b]$ denote the finite set
of integers $x$ with $a \leq x \leq b$. For $n\in\N$ we define the {\em 
(backward) stretch 
of} $A$ with respect to $n$ by $$\str(A,n) = 0 \text{ if } n\notin A,$$ and   
$$\str(A) = \max\{k\leq n: [n-k,n]
\subseteq A\} + 1 \text{ if } n\in A.$$

We first illustrate and motivate graphically the notion of $\str(A,n)$, 
as well as the construction of $A'$, 
for the example $A = \{3, 4, 5, 10, 12\} \in \PN$. 

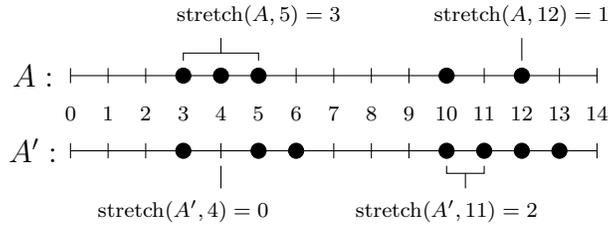
\begin{center}
\begin{tikzpicture}

\draw (0,0) -- (7,0);
\draw (0,-1) -- (7,-1);

\draw (0,0.1) -- (0,-0.1);
\draw (0,-0.9) -- (0,-1.1);
\draw (0.5,0.1) -- (0.5,-0.1);
\draw (0.5,-0.9) -- (0.5,-1.1); 
\draw (1,0.1) -- (1,-0.1);
\draw (1,-0.9) -- (1,-1.1);
\draw (1.5,0.1) -- (1.5,-0.1);
\draw (1.5,-0.9) -- (1.5,-1.1); 
\draw (2,0.1) -- (2,-0.1);
\draw (2,-0.9) -- (2,-1.1);
\draw (2.5,0.1) -- (2.5,-0.1);
\draw (2.5,-0.9) -- (2.5,-1.1); 
\draw (3,0.1) -- (3,-0.1);
\draw (3,-0.9) -- (3,-1.1);
\draw (3.5,0.1) -- (3.5,-0.1);
\draw (3.5,-0.9) -- (3.5,-1.1); 
\draw (4,0.1) -- (4,-0.1);
\draw (4,-0.9) -- (4,-1.1);
\draw (4.5,0.1) -- (4.5,-0.1);
\draw (4.5,-0.9) -- (4.5,-1.1); 
\draw (5,0.1) -- (5,-0.1);
\draw (5,-0.9) -- (5,-1.1);
\draw (5.5,0.1) -- (5.5,-0.1);
\draw (5.5,-0.9) -- (5.5,-1.1);
\draw (6,0.1) -- (6,-0.1);
\draw (6,-0.9) -- (6,-1.1);
\draw (6.5,0.1) -- (6.5,-0.1);
\draw (6.5,-0.9) -- (6.5,-1.1);
\draw (7,0.1) -- (7,-0.1);
\draw (7,-0.9) -- (7,-1.1);


\node at (-0.5, 0) {$A:$};
\node at (-0.5, -1) {$A':$};

\node at (0,-0.5) {\tiny 0};
\node at (0.5,-0.5) {\tiny 1};
\node at (1,-0.5) {\tiny 2};
\node at (1.5,-0.5) {\tiny 3};
\node at (2,-0.5) {\tiny 4};
\node at (2.5,-0.5) {\tiny 5};
\node at (3,-0.5) {\tiny 6};
\node at (3.5,-0.5) {\tiny 7};
\node at (4,-0.5) {\tiny 8};
\node at (4.5,-0.5) {\tiny 9};
\node at (5,-0.5) {\tiny 10};
\node at (5.5,-0.5) {\tiny 11};
\node at (6,-0.5) {\tiny 12};
\node at (6.5,-0.5) {\tiny 13};
\node at (7,-0.5) {\tiny 14};


\filldraw[black] (1.5,0) circle (1mm) node[anchor=west]{};
\filldraw[black] (2,0) circle (1mm) node[anchor=west]{};
\filldraw[black] (2.5,0) circle (1mm) node[anchor=west]{};
\filldraw[black] (5,0) circle (1mm) node[anchor=west]{};
\filldraw[black] (6,0) circle (1mm) node[anchor=west]{};


\filldraw[black] (1.5, -1) circle (1mm) node[anchor=west]{};
\filldraw[black] (2.5, -1) circle (1mm) node[anchor=west]{};
\filldraw[black] (3, -1) circle (1mm) node[anchor=west]{};
\filldraw[black] (5, -1) circle (1mm) node[anchor=west]{};
\filldraw[black] (5.5, -1) circle (1mm) node[anchor=west]{};
\filldraw[black] (6, -1) circle (1mm) node[anchor=west]{};
\filldraw[black] (6.5, -1) circle (1mm) node[anchor=west]{};


\draw (1.5, 0.2) -- (1.5, 0.3) -- (2.5, 0.3) -- (2.5, 0.2);
\draw (2, 0.3) -- (2, 0.5); 

\draw (5, -1.2) -- (5, -1.3) -- (5.5, -1.3) -- (5.5, -1.2);
\draw (5.25, -1.3) -- (5.25, -1.5); 

\draw (6, 0.2) -- (6, 0.5); 

\draw (2, -1.2) -- (2, -1.5); 


\node at (2.5, 0.8) {\tiny $\str(A, 5) = 3$};
\node at (5, -1.8) {\tiny $\str(A', 11) = 2$};
\node at (6, 0.8) {\tiny $\str(A, 12) = 1$};
\node at (1.5, -1.8) {\tiny $\str(A', 4) = 0$};

\end{tikzpicture}

\end{center}

We can quickly verify that 
for $A' = \{3, 5, 6, 10, 11, 12, 13\}$ we have $A\oplus A' = \emptyset$.

\begin{proposition}\label{inverseprop}
{\em Let $A\in\PN$ be non-empty, and let 

$A' = \{x \in A: \str(A, x) \text{ is odd}\} \; \cup $ \\
$ \hspace*{1cm}  \{y \in \N\setminus A: y > 0 \text{ and stretch}(A,y-1) 
\text{ is odd}\}.$

Then $A\oplus A' = \emptyset$.}
\end{proposition}

({\sl Note that we call $n\in\N$} odd {\sl if $n = 2k+1$ for 
some $k\in\N$}.)

{\bf Proof of \ref{inverseprop}.} By fact \ref{basicfact} we need to show
that $$A \, \triangle \, A' = (A\cap A') + 1.$$ In the following we
show that either set is a subset of the other set.

$\subseteq$: Suppose that $x \in A \, \triangle \, A'$.

{\sl Case 1.1} : $x \in A \setminus A'$. This means that
$\str(A,x) > 0$. From $x\notin A'$ and the definition
of $A'$ we get that $\str(A,x)$ is {\em even}. So
in particular $\str(A,x) \geq 2$, implying $x-1 \in A$ 
and $x\geq 1$. Therefore
$\str(A, x-1)$ is {\em odd}, implying $x-1\in A'$ by
definition of $A'$. So $x-1 \in (A\cap A')$, whence $x \in (A\cap A') +1$.

{\sl Case 1.2} : $x \in A' \setminus A$. By definition of $A'$,
this means that $x > 0$ and $\str(A, x-1)$ is odd. So
$x-1 \in A$, and the definition of $A'$ implies $x-1 \in A'$,
yielding $x-1 \in (A\cap A')$ and $x \in (A\cap A')+1$. 

$\supseteq$: Suppose that $x \in (A\cap A') + 1$. In particular,
$x > 0$ and $x-1 \in (A\cap A')$. The statements $x-1 \in A'$ and $x-1\in A$ 
and the definition of $A'$ collectively give us:

\begin{quote}
$(\star)$ \hspace*{4mm} $\str(A,x-1)$ is odd.
\end{quote}

{\sl Case 2.1} : $x\in A$. Statement $(\star)$ and the definition of
$\str(\cdot, \cdot)$ imply $\str(A, x)$ is {\em even}, and by
the definition of $A'$ we get $x \notin A'$. So we get $x\in A\setminus A'$.

{\sl Case 2.2} : $x\notin A$. The definition of $A'$ and $(\star)$
jointly imply $x\in A'$, therefore $x\in A'\setminus A$.

So we established that $A \, \triangle \, A' = (A\cap A') + 1$, which
is equivalent to $A\oplus A' = \emptyset$. \hfill{$\Box$}

\section{Further inquiries}

\subsection{Uniquess of solutions to $A \oplus X = B$} I think that 
the inverses constructed in proposition \ref{inverseprop} are unique, 
and there could be an inductive argument showing this. Moreover, it seems
that the following more general statement holds:
\begin{quote} For all $A, B\in \PN$ there is a unique $X\in \PN$ such that 
$A \oplus X = B$. 
\end{quote}
\subsection{Associative substructures of $(\PN,\oplus)$} One interesting
direction in the analysis of $\oplus$ is the search for ``sub-groups'', 
that is, associative subsets of $\PN$ closed under $\oplus$ and inverses.
Which finite or infinite Abelian groups are isomorphic to
a sub-group of $(\PN,\oplus)$? 
Moreover, Zorn's Lemma implies that every sub-group of $(\PN,\oplus)$ 
is contained in a maximal sub-group with respect to set inclusion
$\subseteq$. Does every maximal subgroup of $(\PN, \oplus)$ 
have the same cardinality?

\end{document}